\documentclass[12pt,reqno]{amsart}

\usepackage{amsmath,amsthm,amsfonts,amscd,amssymb,euscript}
\numberwithin{equation}{section} 
\newcommand{\dbar}{\ensuremath{\bar \partial}}

\newcommand{\C}{\ensuremath{{\mathbb C}}}
\newcommand{\pj}{\ensuremath{{\mathbb P}}}

\newcommand{\N}{\ensuremath{{\mathbb N}}}

\newcommand{\smooth}{\ensuremath{C^{\infty}}}

\newcommand{\I}{\ensuremath{\mathcal I}}

\newcommand{\V}{\ensuremath{\mathcal V}}

\newcommand{\U}{\ensuremath{\mathcal U}}

\newcommand{\oka}{\ensuremath{\mathcal O}}

\newcommand{\cv}{\ensuremath{\mathcal C}}

\begin{document}
\title[On the Relationship between D'Angelo $q$-type and Catlin $q$-type]{On the Relationship between D'Angelo $q$-type and \\ Catlin $q$-type}
\author{Vasile Brinzanescu}

\address{Simion Stoilow Institute of Mathematics of the Romanian Academy, Research unit 3, 21 Calea Grivitei Street, 010702 Bucharest, Romania}

\email{Vasile.Brinzanescu@imar.ro}

\author{Andreea C. Nicoara}

\address{Department of Mathematics, University of Pennsylvania, 209 South $33^{rd}$ St.,  Philadelphia, PA 19104, USA}

\email{anicoara@math.upenn.edu}

\subjclass[2010]{Primary 32F18; 32T25; Secondary 32V35; 13H15.}

\keywords{orders of contact, D'Angelo finite q-type, Catlin finite q-type, finite type domains in $\C^n,$ pseudoconvexity}

\begin{abstract}
We establish inequalities relating two measurements of the order of contact of q-dimensional complex varieties with a real hypersurface.
\end{abstract}

\maketitle

\tableofcontents

\section{Introduction}

The study of the order of contact of complex varieties with the boundary of a domain in $\C^n$ stems from the investigation of the subellipticity of the $\dbar$-Neumann problem. Kohn proved in 1979 in \cite{kohnacta} that for a pseudoconvex domain in $\C^n$ with real-analytic boundary the subellipticity of the $\dbar$-Neumann problem for $(p,q)$ forms is equivalent to the property that all holomorphic varieties of complex dimension $q$ have finite order of contact with the boundary of the domain. D'Angelo introduced a quantitative measure, written $\Delta_q,$ for this order of contact. D'Angelo fleshed out its more important properties culminating with openness and finite determination, which he established in 1982 in \cite{opendangelo}. Meanwhile, Catlin extended Kohn's result to smooth pseudoconvex domains in \cite{catlinnec}, \cite{catlinbdry}, and \cite{catlinsubell}. The notion of finite order of contact of holomorphic varieties of complex dimension $q$ with the boundary of the domain that he defined in \cite{catlinsubell} and showed is equivalent to the subellipticity of the $\dbar$-Neumann problem for $(p,q)$ forms for smooth pseudoconvex domains is not the same as D'Angelo's notion. The two trivially agree for $q=1,$ i.e. for holomorphic curves, but for $q>1$ Catlin merely expressed the hope that they might be shown to equal each other. Catlin's notion is what became known as Catlin $q$-type, $D_q.$ In 1999 in a joint survey paper by D'Angelo and Kohn \cite{kohndangelo}, it was claimed that these two notions ought to be simultaneously finite.

In \cite{catlinsubell} Catlin also proved a lower bound for subelliptic gain in the $\dbar$-Neumann problem $$\epsilon \geq \tau^{-n^2 \, \tau^{n^2}}$$ that holds for any smooth pseudoconvex domain in $\C^n$ and is exponential in $\tau=D_q,$ his notion of contact of holomorphic varieties of complex dimension $q$ with the boundary of the domain. Apart from Catlin's result, there are a number of either sharp or effective bounds for subelliptic gain for $(0,1)$ forms, i.e. when $q=1,$ in terms of $\Delta_1=D_1;$ see \cite{siunote}, \cite{catlincho}, \cite{khanhzampieri}, and \cite{fribourgcatda}. Any other such result for $q>1$ obtained in terms of D'Angelo's more standardly used notion of $q$-type would have to be compared against Catlin's benchmark estimate. Herein lies the significance of our work in this paper as we relate $\Delta_q$ to $D_q$ for $q>1,$ thus enabling this type of comparison.

Both $\Delta_q$ and $D_q$ can also be defined for ideals $\I$ in the ring $\oka_{x_0}$ of germs of holomorphic functions at $x_0.$ Comparing $\Delta_q$ with $D_q$ for such an ideal is much simpler, so we first prove a result of this nature:

\newtheorem{idealthm}{Theorem}[section]
\begin{idealthm}
Let $\I$ be an ideal of germs of holomorphic functions at $x_0,$ then for $1 \leq q \leq n$ $$D_q(\I, x_0) \leq \Delta_q(\I, x_0) \leq \left(D_q(\I, x_0)\right)^{n-q+1}.$$
\label{idealtheorem}
\end{idealthm}

\noindent Theorem~\ref{idealtheorem} is consistent with the simple result that $\Delta_n(\I, x_0)=D_n(\I, x_0)$ for any ideal $\I$ of germs of holomorphic functions in $n$ variables.

\smallskip Our main result is the following:

\newtheorem{mainthm}[idealthm]{Theorem}
\begin{mainthm}
\label{maintheorem}
 Let $\Omega$ in $\C^n$ be a domain with $\smooth$
boundary. Let $x_0 \in b \Omega$ be a point on the boundary of
the domain, and let $1 \leq q <n.$
\begin{enumerate}
\item[(i)] $D_q(b \Omega, x_0) \leq \Delta_q(b \Omega, x_0);$
\item[(ii)] If $\Delta_q(b \Omega, x_0)<\infty$ and the domain is $q$-positive at $x_0$ (the $q$ version of D'Angelo's property P), then $$\Delta_q(b \Omega, x_0) \leq 2 \left(\frac{ D_q(b \Omega, x_0)}{2} \right)^{n-q}.$$
\end{enumerate}
 In particular, if $b \Omega$ is pseudoconvex at $x_0$ and $\Delta_q(b \Omega, x_0)<\infty,$ then $$D_q(b \Omega, x_0) \leq \Delta_q(b \Omega, x_0) \leq 2 \left(\frac{ D_q(b \Omega, x_0)}{2} \right)^{n-q}.$$
\end{mainthm}

\medskip
Since $\Delta_q(b \Omega, x_0)=D_q(b \Omega, x_0)$ for $q=1,$ inequality (i) is sharp.
By definition, $ D_q(b \Omega, x_0) \geq 2,$ and $$2 \left(\frac{ D_q(b \Omega, x_0)}{2} \right)^{n-q} =   D_q(b \Omega, x_0) \,  \left(\frac{ D_q(b \Omega, x_0)}{2} \right)^{n-q-1},$$ so Theorem~\ref{maintheorem} (ii) is not sharp. It is, however, the best result that can be obtained given our method. An example illustrating this point will be provided in Section~\ref{mainthmpf}.  We exclude the value $q=n$ because $b \Omega$ has real dimension $2n-1,$ so looking at its order of contact with an $n$ dimensional complex variety does not make sense. It is also known that subellipticity with exponent $\epsilon =1$ holds at all boundary points for $(p,n)$ forms. The reader may consult p.83 of \cite{kohnacta}. The last part of Theorem ~\ref{maintheorem} follows because a pseudoconvex domain where $\Delta_q(b \Omega, x_0)<\infty$ satisfies $q$-positivity at $x_0,$ a generalization of D'Angelo's property P for $q>1.$ We are deliberately avoiding the terminology property P here in order to be consistent with D'Angelo's usage in \cite{dangelo}. D'Angelo introduced property P in \cite{opendangelo} for a notion of positivity more general than pseudoconvexity. Shortly afterward, Catlin introduced Property (P) in \cite{catlinpropp}, which has since become a standard notion in several complex variables. Details can be found in \cite{catlinpropp}, \cite{boasstraubesurvey}, and \cite{straubebook}. The two names are similar enough to create confusion, so D'Angelo suppressed the term property P in subsequent work, a practice we are following here by employing $q$-positivity instead. 
We would also like to note that our method of proving Theorem~\ref{maintheorem} (ii) breaks down completely in the absence of $q$-positivity, and we have no examples on which we could even formulate a conjecture as to whether $\Delta_q$ are $D_q$ remain simultaneously finite.

The paper is organized as follows: Section~\ref{findap} defines D'Angelo $q$-type and outlines a number of its properties. D'Angelo's property P is also defined here along with $q$-positivity, its $q$ version for $q>1.$ Section~\ref{fincat} is devoted to the Catlin $q$-type. The two notions are then related to each other in Section~\ref{mainthmpf}, where Theorems~\ref{idealtheorem} and \ref{maintheorem} are also proven.

\smallskip
\noindent {\bf Acknowledgements}
The authors wish to thank Catlin and D'Angelo for a number of essential discussions. Additionally, the authors are very grateful to the referee for his suggestions that greatly improved this paper. The first author was partially supported by a grant of the Ministry of National Education, CNCS-UEFISCDI, project number PN-II-ID-PCE-2012-4-0156. He would like to thank the Department of Mathematics at the University of Pennsylvania for the hospitality during the preparation of part of this article.

\section{D'Angelo $q$-type and $q$-positivity}
\label{findap}

Starting with \cite{dangelofirst}, D'Angelo introduced various numerical functions that measure the maximum order of contact of holomorphic varieties of complex dimension $q$ with a real hypersurface $M$ in $\C^n$ such as the boundary of a domain; see \cite{dangelo}.

We shall first give the classical definition of order of contact for $q=1,$ holomorphic curves. Let $r$ be a defining function for the real hypersurface $M$ in $\C^n$. Let $\cv=\cv(m,p)$ be the set of all germs of holomorphic curves $$\varphi: (U,0) \rightarrow (\C^m, p),$$ where $U$ is some neighborhood of the origin in $\C^1$ and $\varphi(0)=p.$ For all $t \in U,$ $\varphi(t)=(\varphi_1(t), \dots, \varphi_m(t)),$ where $\varphi_j(t)$ is holomorphic for every $j$ with $1 \leq j \leq m.$ For each component $\varphi_j,$ the order of vanishing at the origin ${\text ord}_0 \, \varphi_j$ is the order of the first non-vanishing derivative of $\varphi_j,$ i.e. $s \in \N$ such that $$\frac{d}{dt} \varphi_j (0) = \cdots = \frac{d ^{s-1}}{dt^{s-1}} \varphi_j (0) = 0,$$ but $\frac{d ^s}{dt^s} \varphi_j (0) \neq 0.$ We set ${\text ord}_0 \, \varphi = \min_{1 \leq j \leq m} \, {\text ord}_0 \, \varphi_j.$ Consider $\varphi^* r,$ the pullback of $r$ to $\varphi,$ and let $ {\text ord}_0 \, \varphi^* r$ be the order of the first non-vanishing derivative at the origin of $ \varphi^* r$ viewed as a function of $t.$

\medskip
\newtheorem{firstft}{Definition}[section]
\begin{firstft}
Let $M$ be a real hypersurface in $\C^n,$ and let $r$ be a defining function for $M.$ The D'Angelo $1$-type at $x_0 \in M$ is given by $$\Delta_1 (M, x_0) = \sup_{\varphi \in \cv(n,x_0)} \frac{ {\text ord}_0 \, \varphi^* r}{{\text ord}_0 \, \varphi}.$$If $\Delta_1 (M, x_0)$ is finite, we call $x_0$ a point of finite D'Angelo $1$-type.
\end{firstft}

\smallskip When holomorphic varieties have complex dimension greater than $1,$ there is no longer just one natural definition of their order of contact with a real hypersurface in $\C^n$ as not every holomorphic variety of dimension $q \geq 2$ has a local parametrization. Following D'Angelo in \cite{opendangelo}, one approach is to reduce this case to computing $\Delta_1 (\tilde M, x_0)$ for a related hypersurface $\tilde M$ sitting in a different $\C^m$ such that the holomorphic varieties of dimension $q$ generically become holomorphic curves in the new ambient space. Let $\phi : \C^{n-q+1} \rightarrow \C^n$ be any linear embedding of $\C^{n-q+1}$ into $\C^n.$ For generic choices of $\phi,$ the pullback $\phi^* M$ will be a hypersurface in $\C^{n-q+1}.$ We can thus define $\Delta_q (M, x_0)$ as follows:

\medskip
\newtheorem{qft}[firstft]{Definition}
\begin{qft}
\label{qfinitetype}
Let $M$ be a real hypersurface in $\C^n,$ and let $r$ be a defining function for $M.$ The D'Angelo $q$-type at $x_0 \in M$ is given by $$\Delta_q (M, x_0) =\inf_\phi \sup_{\varphi \in \cv(n-q+1,x_0)} \frac{ {\text ord}_0 \, \varphi^* \phi^*r}{{\text ord}_0 \, \varphi}=\inf_\phi \Delta_1 (\phi^*r, x_0),$$ where $\phi : \C^{n-q+1} \rightarrow \C^n$ is any linear embedding of $\C^{n-q+1}$ into $\C^n$ and we have identified $x_0$ with $\phi^{-1} (x_0).$ If $\Delta_q (M, x_0)$ is finite, we call $x_0$ a point of finite D'Angelo $q$-type.
\end{qft}

\medskip
\newtheorem{propdeltaq}[firstft]{Theorem}
\begin{propdeltaq}
Let $M$ be a smooth real hypersurface in $\C^n.$
\begin{enumerate}
\item[(i)] $\Delta_q (M, x_0)$ is well-defined, i.e. independent of the defining function $r$ chosen for $M.$
\item[(ii)] $\Delta_q (M, x_0)$ is not upper semi-continuous in general; see \cite{nonusc}.
\item[(iii)]  Let $\Delta_q (M, x_0)$ be finite at some $x_0 \in M,$ then there exists a neighborhood $V$ of $x_0$ on which $$\Delta_q (M, x) \leq 2 (\Delta_q(M, x_0))^{n-q}.$$
\item[(iv)] The function $\Delta_q(M, x_0)$ is finite determined. In other words, if $\Delta_q(M, x_0)$ is finite, then there exists an integer $k$ such that $\Delta_q(M, x_0)=\Delta_q(M',x_0)$ for $M'$ a hypersurface defined by any $r'$ that has the same $k$-jet at $x_0$ as the defining function $r$ of $M.$
\end{enumerate} \label{propdeltaqthm}
\end{propdeltaq}

\medskip\noindent {\bf Remarks:}

\noindent (1) Part (iii) is Theorem 6.2 from p.634 of \cite{opendangelo} and implies the set of points of finite $q$-type is open. Note that the result holds independently of pseudoconvexity.

\noindent (2) Part (iv) is Proposition 14 from p.88 of \cite{dangeloftc} whose proof implies that if $t = \Delta_q (M, x_0) < \infty,$ then  we can let $k = \lceil t \rceil,$ the ceiling of $t,$ i.e. the least integer greater than or equal to $t.$

\medskip\noindent For the purpose of relating  $\Delta_q(b \Omega, x_0)$ with  $D_q(b \Omega, x_0),$ we will need to show that $\Delta_q(b \Omega, x_0)$ is generic with respect to the choices of linear embeddings $\phi : \C^{n-q+1} \rightarrow \C^n.$ In fact, linear embeddings $\phi : \C^{n-q+1} \rightarrow \C^n$ are in one-to-one correspondence with non-degenerate sets of $q-1$ linear forms $\{w_1, \dots, w_{q-1} \}$ in $\oka_{x_0},$ the local ring of holomorphic germs in $n$ variables at $x_0 \in \C^n.$ The zero set of $\{w_1, \dots, w_{q-1} \}$ is locally the image of the embedding $\phi.$

Restating the embedding $\phi$ as a non-degenerate set of linear forms points to the necessity of having a notion of type that applies to an ideal rather than just a hypersurface, which is what has been defined so far. Indeed, D'Angelo makes the following definition on p.86 of \cite{dangeloftc}:

\medskip
\newtheorem{idtype}[firstft]{Definition}
\begin{idtype}
\label{idealtype}
Let $C^\infty_{x_0}$ be the ring of smooth germs at $x_0 \in \C^n$ and let $\I$ be an ideal in $C^\infty_{x_0}.$ $$\Delta_1(\I,x_0)   = \sup_{\varphi \in \cv(n,x_0)} \:\: \inf_{g \in \I} \:\:\frac{ {\text ord}_0 \, \varphi^* g}{{\text ord}_0 \, \varphi}.$$
\end{idtype}

\medskip\noindent {\bf Remark:} If $M$ be a real hypersurface in $\C^n$ and $x_0 \in M,$ let $\I(M)$ be the ideal of smooth germs in $C^\infty_{x_0}$ that vanish on the germ of $M$ at $x_0.$ Then $\Delta_1(\I(M) ,x_0) = \Delta_1 (M, x_0)$ because the infimum in Definition~\ref{idealtype} is realized by a defining function of $M,$ which has order $1$ at $x_0.$

\medskip Now we can give an equivalent definition to Definition~\ref{qfinitetype} that was first stated by D'Angelo on the bottom of p.86 of \cite{dangeloftc}:

\medskip
\newtheorem{qftequiv}[firstft]{Definition}
\begin{qftequiv}
\label{qfinitetypeequiv}
Let $M$ be a real hypersurface in $\C^n,$ and let $x_0 \in M.$ The D'Angelo $q$-type at $x_0 \in M$ is given by $$\Delta_q (M, x_0) =\inf_{\{w_1, \dots, w_{q-1} \}} \Delta_1\Big( (\I(M), w_1, \dots, w_{q-1}), x_0 \Big),$$ where $\{w_1, \dots, w_{q-1} \}$ is a non-degenerate set of linear forms in $\oka_{x_0},$ $(\I(M), w_1, \dots, w_{q-1})$ is the ideal in $C^\infty_{x_0}$ generated by $\I(M), w_1, \dots, w_{q-1},$ and the infimum is taken over all such non-degenerate sets $\{w_1, \dots, w_{q-1} \}$ of linear forms in $\oka_{x_0}.$
\end{qftequiv}

\medskip This same definition can also be given for an ideal $\I$ in $\oka_{x_0}$ and is the notion that appears in the statement of Theorem~\ref{idealtheorem}: 

\medskip
\newtheorem{qftideal}[firstft]{Definition}
\begin{qftideal}
\label{qfinitetypeideal}
If $\I$ is an ideal in $\oka_{x_0},$
$$\Delta_q(\I, x_0)=\inf_{\{w_1, \dots, w_{q-1} \}} \Delta_1\Big( (\I, w_1, \dots, w_{q-1}), x_0 \Big)=\inf_{\{w_1, \dots, w_{q-1} \}}\sup_{\varphi \in \cv(n,x_0)} \:\: \inf_{g \in (\I, w_1, \dots, w_{q-1})} \:\:\frac{ {\text ord}_0 \, \varphi^* g}{{\text ord}_0 \, \varphi},$$ where $\{w_1, \dots, w_{q-1} \}$ is a non-degenerate set of linear forms in $\oka_{x_0},$ $(\I, w_1, \dots, w_{q-1})$ is the ideal in $\oka_{x_0}$ generated by $\I, w_1, \dots, w_{q-1},$ and the infimum is taken over all such non-degenerate sets $\{w_1, \dots, w_{q-1} \}$ of linear forms in $\oka_{x_0}.$
\end{qftideal}

\medskip Since working in the ring $C^\infty_{x_0}$ is not particularly easy, it would be helpful to reduce the computation of $\Delta_q (M, x_0)$ to a computation in $\oka_{x_0},$ which has much better algebraic properties. Let us assume for the moment that $\Delta_q (M, x_0) = t < \infty,$ and let $k = \lceil t \rceil,$ the ceiling of $t.$ By Theorem~\ref{propdeltaqthm} (iv) and remark (2) following it, $\Delta_q (M, x_0) = \Delta_q (M_k, x_0),$ where $M_k$ is real hypersurface defined by $r_k,$ the polynomial that has the same $k$-jet at $x_0$ as the defining function $r$ of $M.$ The advantage of working with $r_k$ is that we can apply polarization to it, namely we can give a holomorphic decomposition for $r_k$ as $$r_k = Re\{h\} + ||f||^2-||g||^2,$$ where $||f||^2= \sum_{j=1}^N |f_j|^2,$ $||g||^2= \sum_{j=1}^N |g_j|^2,$ and the functions $h, f_1, \dots, f_N, g_1, \dots, g_N$ are all holomorphic polynomials in $n$ variables. This idea first appeared in Section III of \cite{opendangelo}. Furthermore, if $\U(N)$ is the group of $N \times N$ unitary matrices, then for every such unitary matrix $U \in \, \U(N),$ we can consider the ideal of holomorphic polynomials $\I(U, x_0)=(h, f-Ug)$ generated by $h$ and the $N$ components of $f-Ug,$ where $f=(f_1, \dots, f_N)$ and $g= (g_1, \dots, g_N).$ It turns out that
\begin{equation}
\label{delta1bound}
\sup_{U \in \, \U(N)} \Delta_1(\I(U,x_0),x_0) \leq \Delta_1(M_k, x_0) \leq 2 \sup_{U \in \, \U(N)}  \Delta_1(\I(U,x_0),x_0).
\end{equation} 
We have used Corollary 3.7 on p.627 of \cite{opendangelo}. If we now combine this result with Definition~\ref{qfinitetypeequiv} and Theorem~\ref{propdeltaqthm} (iv), we obtain
\begin{equation}
\begin{split}
\label{iuvsdeltaeq}
&\inf_{\{w_1, \dots, w_{q-1} \}} \sup_{U \in \, \U(N)} \Delta_1\Big( (\I(U,x_0), w_1, \dots, w_{q-1}), x_0 \Big) \leq \Delta_q(M, x_0) \\ &\leq 2\inf_{\{w_1, \dots, w_{q-1} \}} \sup_{U \in \, \U(N)} \Delta_1\Big( (\I(U,x_0), w_1, \dots, w_{q-1}), x_0 \Big).
\end{split}
\end{equation}
Since $\U(N)$ is compact, the supremum over it would be easy to handle if $\Delta_1$ were an upper semi-continuous quantity, but it is not as shown by D'Angelo in \cite{nonusc}. We thus need to compare $\Delta_1$ to some other quantity computed in $\oka_{x_0}$ that is upper semi-continuous. We will use $$D(\I, x_0) = \dim_\C (\oka_{x_0} / \I),$$ where $\I$ is an ideal of holomorphic germs at $x_0.$ Here $\dim_\C (\oka_{x_0} / \I)$ means the dimension of $\oka_{x_0} / \I$ viewed as a vector space over $\C.$ This notion appears under different names in the literature. For example, on p.153 of \cite{catlinsubell},  Catlin calls it the multiplicity of the ideal $\I.$

\medskip
\newtheorem{dusc}[firstft]{Proposition}
\begin{dusc}
\label{duscprop}
Let $\I(\lambda)$ be an ideal in $\oka_{x_0}$ that depends continuously on $\lambda.$ Then $D\big(\I(\lambda), x_0\big)$ is an upper semi-continuous function of $\lambda.$ 
\end{dusc}

\noindent This result is part of Proposition 5.3 on p.39 of \cite{tougeron} cited by D'Angelo in \cite{opendangelo}. In our case, $\I(U, x_0)$ obviously depends continuously on $U,$ so $D\big(\I(U,x_0), x_0\big)$ is upper semi-continuous on the compact set $\U(N).$ Thus $D\big(\I(U,x_0), x_0\big)$ achieves a maximum on $\U(N)$ because each $U \in \, \U(N)$ has an open neighborhood $V(U)$ such that $$D\big(\I(U',x_0), x_0\big)\leq D\big(\I(U,x_0), x_0\big)$$ for every $U' \in V(U)$ from the upper semi-continuity, $\{V(U)\}$ is an open cover of $\U(N),$ we can thus pass to a finite open subcover $\{V(U_j)\}_{1\leq j \leq p},$ and then we can take $\displaystyle \max_{1\leq j \leq p} D\big(\I(U_j,x_0), x_0\big).$

Since we are primarily interested in the case $q>1,$ the object that appears naturally corresponding to a proper ideal $\I$ in $\oka_{x_0}$ is $D\Big( (\I, w_1, \dots, w_{q-1}), x_0 \Big),$ where $\{w_1, \dots, w_{q-1} \}$ is a non-degenerate set of linear forms. It turns out that $D\Big( (\I, w_1, \dots, w_{q-1}), x_0 \Big)$ is generic when we consider its value over all non-degenerate sets of linear forms $\{w_1, \dots, w_{q-1} \}:$

\medskip
\newtheorem{multgen}[firstft]{Proposition}
\begin{multgen}
\label{multgeneric}
Let $\I$ be a proper ideal in $\oka_{x_0},$ and let $x_0 \in \C^n.$
\begin{equation*}
\inf_{\{w_1, \dots, w_{q-1} \}} D\Big( (\I, w_1, \dots, w_{q-1}), x_0 \Big)={\text gen.val}_{\{w_1, \dots, w_{q-1} \}} D\Big( (\I, w_1, \dots, w_{q-1}), x_0 \Big),
\end{equation*}
where $\{w_1, \dots, w_{q-1} \}$ is a non-degenerate set of linear forms in $\oka_{x_0},$ $(\I, w_1, \dots, w_{q-1})$ is the ideal in $\oka_{x_0}$ generated by $\I, w_1, \dots, w_{q-1},$ and the infimum and the generic value are both taken over all such non-degenerate sets $\{w_1, \dots, w_{q-1} \}$ of linear forms in $\oka_{x_0}.$ In other words, the infimum is achieved and equals the generic value.
\end{multgen}

\smallskip\noindent {\bf Proof:} For every non-degenerate set of linear forms $\{w_1, \dots, w_{q-1} \}$ in $\oka_{x_0},$ we consider a linear change of variables at $x_0$ such that $w_1, \dots, w_{q-1}$ become the coordinate functions $z_1, \dots, z_{q-1}.$ Let $\tilde \I$ be the image of $\I$ under this linear change of variables. Consider now $$D\Big( (\tilde \I, z_1, \dots, z_{q-1}), x_0 \Big) = \dim_\C \Big(\oka_{x_0} / (\tilde \I, z_1, \dots, z_{q-1})\Big).$$ Since the variables $z_1, \dots, z_{q-1}$ get set to zero in the quotient $\oka_{x_0} / (\tilde \I, z_1, \dots, z_{q-1})$ and the quantity $D\Big( (\tilde \I, z_1, \dots, z_{q-1}), x_0 \Big)$ is invariant under linear changes of variables, it follows that $$\inf_{\{w_1, \dots, w_{q-1} \}} D\Big( (\I, w_1, \dots, w_{q-1}), x_0 \Big)={\text gen.val}_{\{w_1, \dots, w_{q-1} \}} D\Big( (\I, w_1, \dots, w_{q-1}), x_0 \Big)$$ as needed. \qed

\medskip Propositions ~\ref{duscprop} and ~\ref{multgeneric} together make looking at $$  \sup_{U \in \, \U(N)} \:\: \inf_{\{w_1, \dots, w_{q-1} \}} \:\: D\Big( (\I(U,x_0), w_1, \dots, w_{q-1}), x_0 \Big)$$ much easier in the sense that the supremum and the infimum, which are both achieved, can be exchanged here. D'Angelo has already related this quantity to $\Delta_q(M, x_0)$ in Theorem 14 on p.91 of \cite{dangeloftc}:

\medskip
\newtheorem{dangeloqbound}[firstft]{Theorem}
\begin{dangeloqbound}
\label{dangeloqboundthm}
\begin{equation*}
\Delta_q(M, x_0) \leq 2   \sup_{U \in \, \U(N)} \:\: \inf_{\{w_1, \dots, w_{q-1} \}} \:\: D\Big( (\I(U,x_0), w_1, \dots, w_{q-1}), x_0 \Big) \leq 2 \Big(\Delta_q(M, x_0) \Big)^{n-q}
\end{equation*}
\end{dangeloqbound}

\noindent We shall use the lower bound, 
\begin{equation}
\label{qbound}
\frac{1}{2} \:\: \Delta_q(M, x_0) \leq  \sup_{U \in \, \U(N)} \:\: \inf_{\{w_1, \dots, w_{q-1} \}} \:\: D\Big( (\I(U,x_0), w_1, \dots, w_{q-1}), x_0 \Big)
\end{equation}
in the proof of Theorem~\ref{maintheorem}.

Another ingredient in the proof of Theorem~\ref{maintheorem} is showing that $\Delta_q(M, x_0)$ assumes the generic value with respect to choices of non-degenerate sets of linear forms $\{w_1, \dots, w_{q-1} \}$ in its definition:

\medskip
\newtheorem{deltaqgen}[firstft]{Proposition}
\begin{deltaqgen}
\label{deltaqgeneric}
Let $M$ be a real hypersurface in $\C^n,$ and let $x_0 \in M.$ The infimum in the definition of $\Delta_q(M, x_0)$ is achieved and equal to the generic value,
\begin{equation*}
\begin{split}
\Delta_q (M, x_0) &=\inf_{\{w_1, \dots, w_{q-1} \}} \:\: \Delta_1\Big( (\I(M), w_1, \dots, w_{q-1}), x_0 \Big)\\&={\text gen.val}_{\{w_1, \dots, w_{q-1} \}} \:\: \Delta_1\Big( (\I(M), w_1, \dots, w_{q-1}), x_0 \Big),
 \end{split}
\end{equation*}
where $\{w_1, \dots, w_{q-1} \}$ is a non-degenerate set of linear forms in $\oka_{x_0},$ $(\I(M), w_1, \dots, w_{q-1})$ is the ideal in $C^\infty_{x_0}$ generated by $\I(M), w_1, \dots, w_{q-1},$ and the infimum and the generic value are both taken over all such non-degenerate sets $\{w_1, \dots, w_{q-1} \}$ of linear forms in $\oka_{x_0}.$
\end{deltaqgen}

\smallskip\noindent {\bf Proof:} We consider two cases. First, if $$\Delta_q (M, x_0) =\inf_{\{w_1, \dots, w_{q-1} \}} \:\: \Delta_1\Big( (\I(M), w_1, \dots, w_{q-1}), x_0 \Big)=\infty,$$ then clearly $\Delta_1\Big( (\I(M), w_1, \dots, w_{q-1}), x_0 \Big)=\infty$ for every non-degenerate set of linear forms $\{w_1, \dots, w_{q-1} \}$ in $\oka_{x_0},$ so $$\Delta_q (M, x_0) ={\text gen.val}_{\{w_1, \dots, w_{q-1} \}} \:\: \Delta_1\Big( (\I(M), w_1, \dots, w_{q-1}), x_0 \Big),$$ i.e. the infimum is achieved and equal to the generic value. Second, assume $\Delta_q (M, x_0) = t <\infty.$ One of the equivalent ways of defining $\Delta_q (M, x_0)$ is as $$\Delta_q (M, x_0) =\inf_{\{w_1, \dots, w_{q-1} \}} \:\:  \sup_{\varphi \in \cv(n,x_0)}  \:\:  \inf_{g \in (\I(M), w_1, \dots, w_{q-1})}  \:\:  \frac{ {\text ord}_0 \, \varphi^* g}{{\text ord}_0 \, \varphi}.$$ Obviously, the supremum is realized here by curves that set $w_1, \dots, w_{q-1}$ to zero as these functions have the lowest order of vanishing, namely $1.$ It follows that $$\Delta_q (M, x_0) ={\text gen.val}_{\{w_1, \dots, w_{q-1} \}} \:\: \Delta_1\Big( (\I(M), w_1, \dots, w_{q-1}), x_0 \Big),$$ and the infimum over all non-degenerate sets $\{w_1, \dots, w_{q-1} \}$ of linear forms in $\oka_{x_0}$ is indeed achieved and equal to the generic value. \qed

\medskip We shall now state as a corollary the equivalent result for $\Delta_q(\I,x_0),$ which is needed for the proof of Theorem~\ref{idealtheorem}:

\medskip
\newtheorem{deltaqgenid}[firstft]{Corollary}
\begin{deltaqgenid}
\label{deltaqgenericideal}
If $\I$ is any ideal in $\oka_{x_0},$ the infimum in the definition of $\Delta_q(\I, x_0)$ is achieved and equal to the generic value,
\begin{equation*}
\begin{split}
\Delta_q (\I, x_0) &=\inf_{\{w_1, \dots, w_{q-1} \}} \:\: \Delta_1\Big( (\I, w_1, \dots, w_{q-1}), x_0 \Big)\\&={\text gen.val}_{\{w_1, \dots, w_{q-1} \}} \:\: \Delta_1\Big( (\I, w_1, \dots, w_{q-1}), x_0 \Big).
 \end{split}
\end{equation*}
\end{deltaqgenid}

\smallskip\noindent {\bf Proof:} The proof of Proposition~\ref{deltaqgeneric} applies verbatim with $\Delta_q (\I, x_0)$ replacing $\Delta_q (M, x_0)$ and $\I$ replacing $\I(M).$ \qed

\medskip For the proof of Theorem~\ref{idealtheorem}, we also need to state part of D'Angelo's Theorem 2.7 from p.622 of \cite{opendangelo}:

\medskip
\newtheorem{dangeloineq}[firstft]{Theorem}
\begin{dangeloineq}
\label{dangeloinequality} Let $\I$ be an ideal of $\oka_{x_0},$ then $\Delta_1(\I,x_0) \leq D(\I,x_0).$
\end{dangeloineq}

\medskip Finally, we turn our attention to D'Angelo's property P, which he defined on p.631 of \cite{opendangelo}:

\medskip
\newtheorem{propP}[firstft]{Definition}
\begin{propP}
\label{propertyP}
Let $M$ be a real hypersurface of $\C^n,$ and let $x_0$ be a point of finite type on $M.$ We suppose that $\Delta_1(M, x_0) <k.$ Let $j_{k,x_0} r = r_k = Re\{h\} + ||f||^2-||g||^2$ be a holomorphic decomposition at $x_0$ of the $k$-jet of the defining function $r$ of $M.$ We say that $M$ satisfies property P at $x_0$ if for every holomorphic curve $\varphi \in \cv(n,x_0)$ for which $\varphi^* h$ vanishes, the following two conditions are satisfied:
\begin{enumerate}
\item[(i)] ${\text ord}_0 \, \varphi^* r$ is even, i.e. ${\text ord}_0 \, \varphi^* r = 2a,$ for some $a \in \N;$
\item[(ii)] $\displaystyle \left(\frac{d}{dt}\right)^a \left(\frac{d}{d \bar t}\right)^a \varphi^* r (0) \neq 0.$
\end{enumerate}
\end{propP}

\medskip\noindent {\bf Remarks:}  

\noindent (1) Due to Theorem~\ref{propdeltaqthm} (iv), the finite determination property of $\Delta_q(M, x_0)$ for all $1 \leq q <n,$ this definition is independent of $k$ provided $k$ is large enough. 

\noindent (2) Since the D'Angelo type does not depend on the coordinate system and is always greater than or equal to $2,$ no holomorphic curve with $\varphi^* h \not\equiv 0$ can realize the supremum. The function $h$ has order $1$ at $x_0$ and can be mapped to $z_n$ via a change of coordinates. Hence holomorphic curves satisfying $\varphi^* h \not\equiv 0$ are irrelevant to the type consideration.

\bigskip\noindent Let us now define the $q$ version of D'Angelo's property P, $q$-positivity, the hypothesis that appears in Theorem~\ref{maintheorem}.

\medskip
\newtheorem{qpropP}[firstft]{Definition}
\begin{qpropP}
\label{qpropertyP}
Let $M$ be a real hypersurface of $\C^n,$ and let $x_0 \in M$ be such that $\Delta_q(M, x_0) <k.$ Let $j_{k,x_0} r = r_k = Re\{h\} + ||f||^2-||g||^2$ be a holomorphic decomposition at $x_0$ of the $k$-jet of the defining function $r$ of $M.$ We say that $M$ is $q$-positive at $x_0$ if for every holomorphic curve $\varphi \in \cv(n,x_0)$ for which $\varphi^* h$ vanishes and such that the image of $\varphi$ locally lies in the zero locus of a non-degenerate set of linear forms $\{w_1, \dots, w_{q-1} \}$ at $x_0,$ the following two conditions are satisfied:
\begin{enumerate}
\item[(i)] ${\text ord}_0 \, \varphi^* r$ is even, i.e. ${\text ord}_0 \, \varphi^* r = 2a,$ for some $a \in \N;$
\item[(ii)] $\displaystyle \left(\frac{d}{dt}\right)^a \left(\frac{d}{d \bar t}\right)^a \varphi^* r (0) \neq 0.$
\end{enumerate}
\end{qpropP}

\medskip\noindent {\bf Remark:} The change here versus D'Angelo's property P is that we ask that his conditions be satisfied only for the holomorphic curves that come into the computation of $\Delta_q(M, x_0).$

\medskip\noindent The reason D'Angelo introduced property P is that it allowed him to prove the following result, which appears as Theorem 5.3 on p.631 of \cite{opendangelo}:

\medskip
\newtheorem{1propPthm}[firstft]{Theorem}
\begin{1propPthm}
\label{1propertyP}
Suppose that $M$ satisfies property P at $x_0,$ then $$\Delta_1(M, x_0) = 2 \,  \Delta_1(\I(U,x_0),x_0),$$ i.e. the upper bound in Equation~\eqref{delta1bound} is achieved.
\end{1propPthm}

A pseudoconvex domain of finite D'Angelo type has property P as do the hypersurfaces corresponding to truncations of the defining function at $x_0$ of any order higher than the type. Before formally stating this result, we shall state a proposition that appeared as Proposition 2 on p.138 of \cite{dangelo}, which justifies why such a result ought to be true.

\medskip
\newtheorem{pscproppaux}[firstft]{Proposition}
\begin{pscproppaux}
\label{pscproppauxiliary}
Suppose that $M$ is a pseudoconvex hypersurface containing the origin with local defining function $r.$ Suppose further that $\varphi: (\C,0) \rightarrow (\C^n, 0)$ is a parametrized holomorphic curve such that the Taylor series for $\varphi^* r$ satisfies
\begin{enumerate}
\item[(i)] ${\text ord}_0 \, \varphi^* r=m$
\item[(ii)] $\displaystyle \left(\frac{d}{dt}\right)^a  \varphi^* r (0) = 0 \quad a \leq m.$
\end{enumerate}
Then the order of vanishing $m=2k$ is even, and the coefficient of $|t|^{2k}$ in $\varphi^* r$ is positive.
\end{pscproppaux}

\medskip\noindent {\bf Remark:} Condition (ii) is eliminating pure terms up to and including of order $m,$ which is the vanishing order of $ \varphi^* r.$ The same is achieved via the requirement that  $\varphi^* h \equiv 0$ in the definition of D'Angelo's property P as well as in the definition of $q$-positivity.

\medskip\noindent The following result appears on p.632 of \cite{opendangelo}:

\medskip
\newtheorem{pscpropp}[firstft]{Proposition}
\begin{pscpropp}
\label{pscproppproposition}
Suppose $M$ is pseudoconvex near $x_0,$ and that $\Delta_1(M, x_0)$ is finite. Then $M$ and $M_k,$ the hypersurface corresponding to the truncation of order $k$ of the defining function at $x_0,$ satisfy property P at $x_0$ for all sufficiently large $k.$
\end{pscpropp}

It is easy to see from Proposition~\ref{pscproppauxiliary} that the equivalent result should hold for $q$-positivity as well:

\medskip
\newtheorem{pscqpropp}[firstft]{Proposition}
\begin{pscqpropp}
\label{pscqproppproposition}
Suppose $M$ is pseudoconvex near $x_0,$ and that $\Delta_q (M, x_0)$ is finite. Then $M$ and $M_k,$ the hypersurface corresponding to the truncation of order $k$ of the defining function at $x_0,$ are $q$-positive at $x_0$ for all sufficiently large $k.$
\end{pscqpropp}

\section{Catlin $q$-type}
\label{fincat}

Catlin wished to avoid having to characterize the order of contact of a holomorphic variety $V^q$ of complex dimension $q$ with the boundary of the domain along the singular locus of the variety, which can be considerably more complicated when $q>1$ than for holomorphic curves. To that end, he introduced in \cite{catlinsubell} a numerical function $D_q (M, x_0)$ that measures the order of contact of varieties $V^q$ with $M$ only along generic directions. Following D'Angelo in \cite{opendangelo}, he also defined such an order of contact $D_q(\I, x_0)$ for an ideal $\I$ of holomorphic germs in $\oka_{x_0}.$ 

Let $V^q$ be the germ of a holomorphic variety of complex dimension $q$ passing through $x_0.$ Let $G^{n-q+1}$ be the set of all $(n-q+1)$-dimensional complex planes through $x_0.$ Consider the intersection $V^q \cap S$ for $S \in G^{n-q+1}.$ For a generic, thus open and dense, subset $\tilde W$ of $G^{n-q+1},$ $V^q \cap S$ consists of finitely many irreducible one-dimensional components $V^q_{S,k}$ for $k=1, \dots, P.$ Let us parametrize each such germ of a curve by some open set $U_k \ni 0$ in $\C.$ Thus, $\gamma_S^k : U_k \rightarrow V^q_{S,k},$ where $\gamma_S^k (0)=x_0.$ For every holomorphic germ $f \in \oka_{x_0},$ consider the quantity $$\tau (f, V^q \cap S) = \max_{k=1, \dots, P} \frac{ {\text ord}_0 \, {\left(\gamma^k_S\right)}^* f}{{\text ord}_0 \, \gamma^k_S}.$$ Likewise, for $r$ the defining function of a real hypersurface $M$ in $\C^n$ passing through $x_0,$ set $$\tau (V^q \cap S, x_0) = \max_{k=1, \dots, P} \frac{ {\text ord}_0 \, {\left(\gamma^k_S\right)}^* r}{{\text ord}_0 \, \gamma^k_S}.$$ In Section 3 of \cite{catlinsubell}, Catlin showed $\tau (f, V^q \cap S)$ assumes the same value for all $S$ in a generic subset $\tilde W$ of planes. Therefore, he defined $$\tau(f, V^q) = {\text gen.val}_{S \in \tilde W} \left\{ \tau (f, V^q \cap S) \right\}$$ and $$\tau(\I, V^q) = \min_{f \in \I} \tau(f, V^q).$$

\medskip
\newtheorem{catlinid}{Definition}[section]
\begin{catlinid}
\label{catlinideal}
Let $\I$ be an ideal of holomorphic germs at $x_0,$ then the Catlin $q$-type of the ideal $\I$ is given by $$D_q (\I, x_0)=\sup_{V^q}  \left\{ \tau(\I, V^q) \right\},$$ where the supremum is taken over the set of all germs of $q$-dimensional holomorphic varieties $V^q$ passing through $x_0.$
\end{catlinid}

In the same section 3 of \cite{catlinsubell}, Catlin showed $\tau (V^q \cap S, x_0)$ assumes the same value for all $S$ in a generic subset $\tilde W$ of planes, so he defined $$\tau(V^q, x_0) = {\text gen.val}_{S \in \tilde W} \left\{ \tau (V^q \cap S, x_0) \right\}.$$

\medskip
\newtheorem{catlinft}[catlinid]{Definition}
\begin{catlinft}
\label{catlinfinitetype}
Let $M$ be a real hypersurface in $\C^n.$ The Catlin $q$-type at $x_0 \in M$ is given by $$D_q (M, x_0)=\sup_{V^q}  \left\{ \tau (V^q, x_0) \right\},$$ where the supremum is taken over the set of all germs of $q$-dimensional holomorphic varieties $V^q$ passing through $x_0.$
\end{catlinft}

\noindent Clearly, $\Delta_1 (M, x_0) = D_1 (M , x_0)$ as there is only one $n$-dimensional complex plane passing through $x_0$ in $\C^n.$

Some more explanations are in order regarding Catlin's construction. We have claimed $\tau (f, V^q \cap S)$ and $\tau (V^q \cap S, x_0)$ are constant on a generic set $\tilde W$ of $(n-q+1)$-dimensional complex planes $S$ through $x_0.$ The quantities $\tau (f, V^q \cap S)$ and $\tau (V^q \cap S, x_0)$ are computed by looking at the normalized vanishing orders of $f$  and $r$ respectively along the curves $V^q_{S,k}$ for $k=1, \dots, P$ that represent the intersection of $V^q$ with $S.$ In fact, the number of curves in the intersection, $P$ is the same for all $S \in \tilde W,$ and furthermore, the curves $V^q_{S_a,k}$ can be smoothly parametrized via a parameter $a=(a_1, \dots, a_N)$ for $N=(n-q+1)(q-1),$ the dimension of $\tilde W.$ Proposition 3.1 (ii) on p.140 of \cite{catlinsubell} states that as $S_a \in \tilde W$ varies smoothly, the intersection curves $V^q_{S_a,k}$ do as well, and their number stays constant.

In his proof of Proposition 3.1 from \cite{catlinsubell}, Catlin has to remove three different sets $W_1,$ $W_2,$ and $W_3$ from $G^{n-q+1}$ in order to arrive at his generic set $\tilde W$ on which such strong conclusions hold. To the germ of the variety $V^q,$ there corresponds a prime ideal $\I$ in the ring $\oka_{x_0}$ of all germs of holomorphic functions that vanish on $V^q.$ Catlin uses Gunning's Local Parametrization Theorem from p.16 of \cite{gunning} in order to construct a set of canonical equations for $V^q.$ This construction involves choosing a special set of coordinates where the generators of the ideal simultaneously satisfy the Weierstrass Preparation Theorem with respect to the variables that give the regular system of parameters $\I$ has as a prime ideal in the regular local ring $\oka_{x_0}.$ The intersection $V^q \cap S$ is ill-behaved where $V^q$ does not have pure dimension $q$ as the intersection might consist of points rather than curves as well as along the singular locus of $V^q.$ To remove both, Catlin constructs a conic variety $X'$ whose defining equation consists of the product of the discriminants of the Weierstrass polynomials that give the canonical equations for $V^q$ (these discriminants capture the singular locus of $V^q$) with the additional generator that gives the non pure dimensional part of $V^q.$ $W_1$ consists of all $(n-q+1)$-dimensional complex planes that intersect $X'.$ 

Additionally, for the intersection $V^q \cap S$ to behave well, a good notion of transversality has to apply. Transversality cannot be tested well for curves, which is what $V^q \cap S$ yields generically, but it can be tested very well for points. To reduce the intersection to points, Catlin looks at the conic variety corresponding to $V^q,$ which he calls $V'.$ The variety $V'$ captures the tangent cone of $V^q,$ exactly where singularities of $V^q$ manifest themselves as the dimension of the tangent cone jumps at a singular point. $V'$ still has dimension $q.$ Consider $\tilde V,$ the projective variety in $\pj^{n-1}$ corresponding to $V'.$ $\tilde V$ has dimension $q-1.$ For every $S \in G^{n-q+1},$ there corresponds a projective plane $\tilde S$ of dimension $n-q$ in $\pj^{n-1}.$ Generically, $\tilde V \cap \tilde S$ consists of finitely many points $\tilde z^1, \dots, \tilde z^D$ with transverse intersections, meaning that each $\tilde z^i$ is a smooth point of $\tilde V$ and the tangent spaces satisfy $T_{\tilde z^i} \tilde V \cap T_{\tilde z^i} \tilde S = 0$ for $i=1, \dots, D.$ Let $W_2$ be the subset of $G^{n-q+1}$ where this generic behavior does not take place.

Finally, the construction of $W_1$ involved the use of canonical equations for $V^q.$ Hence variables $z_{q+1}, \dots, z_n$ give the regular system of parameters corresponding to the pure $q$-dimensional part of the variety $V^q.$ The variable $z_q$ corresponds to the additional generator that gives the non pure dimensional part of $V^q.$ The $(n-q+1)$-dimensional complex plane $S$ is defined by the linear equations $ \sum^n_{j=1} \: a^i_j z_j=0$ for $i=1, \dots, q-1,$ which need to be linearly independent. A $(q-1) \times (q-1)$ minor of $(a^i_j)$ should thus have full rank. On the other hand, for the intersection $V^q \cap S$ to behave well, this $(q-1) \times (q-1)$ minor should be exactly $(a^i_j)_{1 \leq i,j \leq q-1}$ with respect to the complementary variables $z_1, \dots, z_{q-1}.$ Therefore, Catlin sets $$W_3= \left\{ S \in G^{n-q+1} \:\: \Big| \:\: \det (a^i_j)_{1 \leq i,j \leq q-1}=0 \right\}.$$

We shall now state the rest of the results from \cite{catlinsubell} that play a role in the proofs of Theorems~\ref{idealtheorem} and \ref{maintheorem}. The following is Catlin's Theorem 3.7 on p.154 of \cite{catlinsubell}:

\medskip
\newtheorem{catlinmult}[catlinid]{Theorem}
\begin{catlinmult}
\label{catlinmultthm}
Let $\I$ be an ideal in $\oka_{x_0},$ then $${\text gen.val}_{\{w_1, \dots, w_{q-1} \}} \:\: D\Big( (\I, w_1, \dots, w_{q-1}), x_0 \Big) \leq \prod_{i=q}^n \: D_i (\I, x_0),$$ where the generic value is computed over all non-degenerate sets $\{w_1, \dots, w_{q-1} \}$ of linear forms in $\oka_{x_0}.$
\end{catlinmult}

\medskip\noindent {\bf Remark:} In the context of Catlin's definitions, the $(n-q+1)$-dimensional complex plane $S$ through $x_0$ is precisely the zero locus of the non-degenerate set of linear forms $\{w_1, \dots, w_{q-1} \}.$

\medskip\noindent Since $D_k(\I, x_0) \leq D_q (\I, x_0)$ for all $k \geq q,$ we obtain the following corollary to Theorem~\ref{catlinmultthm}:

\medskip
\newtheorem{catlinmultcor}[catlinid]{Corollary}
\begin{catlinmultcor}
\label{catlinmultcorollary}
Let $\I$ be an ideal in $\oka_{x_0},$ then $${\text gen. val}_{\{w_1, \dots, w_{q-1} \}}\:\: D\Big( (\I, w_1, \dots, w_{q-1}), x_0 \Big) \leq \left( D_q (\I, x_0)\right)^{n-q+1}.$$
\end{catlinmultcor}

\medskip\noindent In case $\Delta_q (M, x_0)=t<\infty,$ the truncation $r_k$ of the defining function $r$ of $M$ at $x_0$ of order $k = \lceil t \rceil$ has the holomorphic decomposition $r_k = Re\{h\} + ||f||^2-||g||^2,$ and $M$ is $q$-positive at $x_0,$ we would like to relate $\tau(V^q, x_0)$ to $\tau(\I(U), V^q)$ for any unitary matrix $U \in \U(N)$ and any $q$-dimensional complex variety $V^q.$ The reader may wish to check by going through the proof of D'Angelo's Theorem 5.3 on p.631 of \cite{opendangelo}, which we stated here as Theorem~\ref{1propertyP}, that $q$-positivity at $x_0$ implies $$\tau(V^q, x_0) \geq 2 \, \tau(\I(U), V^q) \quad \forall \, U, \: \forall \, V^q.$$ Catlin uses this reasoning on p.156 of \cite{catlinsubell} in order to finish the proof of his Theorem 3.4 without formally defining $q$-positivity. Instead, he employs this argument for a pseudoconvex domain, where we have shown that $q$-positivity holds if the D'Angelo $q$-type is finite. Catlin assumes that $D_q (M, x_0)$ is finite instead.

\section{Comparing $\Delta_q$ with $D_q$}
\label{mainthmpf}

The proof of Theorem~\ref{idealtheorem} is comprised of two results.

\medskip
\newtheorem{idealleft}{Proposition}[section]
\begin{idealleft}
\label{idealleftprop}
Let $\I$ be any ideal in $\oka_{x_0}.$ For any $1 \leq q \leq n,$ $D_q(\I, x_0) \leq \Delta_q(\I, x_0).$
\end{idealleft}

\smallskip\noindent {\bf Proof:} Let $\Delta_q(\I, x_0)=t < \infty,$ else the estimate is trivially true. Assume $D_q(\I, x_0) >t.$ Since $D_q (\I, x_0)$ is defined as the supremum over all $q$-dimensional holomorphic varieties passing through $x_0$ of $\tau (\I,V^q),$ there exists such a holomorphic variety $V^q$ for which $$\tau (\I,V^q)= \min_{f \in \I} \tau(f, V^q) = \min_{f \in \I} {\text gen.val}_{S \in \tilde W} \left\{ \tau (f, V^q \cap S) \right\} =t' >t,$$ but as shown in Proposition~\ref{deltaqgenericideal}, $\Delta_q(\I, x_0)$ is generic over the choice of $S,$ so the curves in $V^q \cap S$ already enter into the computation of $\Delta_q(\I, x_0).$ Therefore, $\Delta_q(\I, x_0) \geq t' > t,$ a contradiction. \qed

\medskip
\newtheorem{idealright}[idealleft]{Proposition}
\begin{idealright}
\label{idealrightprop}
Let $\I$ be any ideal in $\oka_{x_0}.$ For any $1 \leq q \leq n,$ $$\Delta_q(\I, x_0) \leq \left(D_q(\I, x_0)\right)^{n-q+1}.$$
\end{idealright}

\smallskip\noindent {\bf Proof:} Let $\{w_1, \dots, w_{q-1} \}$ be any non-degenerate set of linear forms in $\oka_{x_0}.$ First, we apply Theorem~\ref{dangeloinequality} to the ideal $(\I, w_1, \dots, w_{q-1})$ to obtain $$\Delta_1\Big((\I, w_1, \dots, w_{q-1}),x_0\Big) \leq D\Big((\I, w_1, \dots, w_{q-1}),x_0\Big).$$ Next, we take the generic value over all non-degenerate sets $\{w_1, \dots, w_{q-1} \}$ of linear forms in $\oka_{x_0}:$ $${\text gen.val}_{\{w_1, \dots, w_{q-1} \}}\:\: \Delta_1\Big( (\I, w_1, \dots, w_{q-1}), x_0 \Big) \leq {\text gen.val}_{\{w_1, \dots, w_{q-1} \}} \:\: D\Big((\I, w_1, \dots, w_{q-1}),x_0\Big).$$ By Corollary~\ref{deltaqgenericideal}, $$\Delta_q(\I, x_0) = {\text gen.val}_{\{w_1, \dots, w_{q-1} \}}\:\:  \Delta_1 \Big( (\I, w_1, \dots, w_{q-1}), x_0 \Big),$$ while by Corollary~\ref{catlinmultcorollary}, $${\text gen.val}_{\{w_1, \dots, w_{q-1} \}} \:\: D\Big((\I, w_1, \dots, w_{q-1}),x_0\Big)\leq \left( D_q (\I, x_0)\right)^{n-q+1}.$$
\qed

\medskip\noindent {\bf Proof of Theorem~\ref{idealtheorem}:} Proposition~\ref{idealleftprop} proves the left-hand side inequality, while Proposition~\ref{idealrightprop} proves the right-hand side one. \qed

\medskip As mentioned in the introduction, if $q=n,$ Theorem~\ref{idealtheorem} implies $\Delta_n(\I,x_0)=D_n(\I,x_0)$ for every ideal $\I$ in $\oka_{x_0}.$ From the proof of Proposition~\ref{idealrightprop}, it follows that
\begin{equation}
\label{nval}
\Delta_n(\I, x_0)=D_n(\I,x_0) = {\text gen.val}_{\{w_1, \dots, w_{n-1} \}}\:\:  D \Big( (\I, w_1, \dots, w_{n-1}), x_0 \Big).
\end{equation} 

\medskip\noindent Applying Catlin's Theorem~\ref{catlinmultthm} together with Equation~\eqref{nval} allows us to give the corresponding result to Corollary~\ref{catlinmultcorollary} for an ideal $\I(U)=(h,f-Ug)$ in $\oka_{x_0}$ arising from a holomorphic decomposition $r_k = Re\{h\} + ||f||^2-||g||^2$ of the truncation  $r_k$ of the defining function $r$ of the real hypersurface $M$ at $x_0$ when $\Delta_q (M, x_0)=t$ and $k = \lceil t \rceil:$

\medskip
\newtheorem{catlinmultcorbdry}[idealleft]{Corollary}
\begin{catlinmultcorbdry}
\label{catlinmultbdrycorollary}
Let $M$ be a real hypersurface in $\C^n,$ let $r$ be a defining function for $M,$ and let $x_0 \in M.$ If $\Delta_q (M, x_0)=t,$ $r_k = Re\{h\} + ||f||^2-||g||^2$ is the holomorphic decomposition of the truncation  $r_k$ of the defining function $r$ for $k = \lceil t \rceil,$ and $\I(U)=(h,f-Ug)$ is an ideal in $\oka_{x_0}$ corresponding to this holomorphic decomposition for $U$ a unitary matrix, then $${\text gen. val}_{\{w_1, \dots, w_{q-1} \}}\:\: D\Big( (\I(U), w_1, \dots, w_{q-1}), x_0 \Big) \leq \left( D_q (\I(U), x_0)\right)^{n-q}.$$
\end{catlinmultcorbdry}

\smallskip\noindent {\bf Proof:} $r_k = Re\{h\} + ||f||^2-||g||^2$ defines a real hypersurface in $\C^n,$ so the holomorphic function $h$ must contain a term of first order, else the gradient of $r_k$ would be zero at $x_0.$ Therefore, the ideal $\I(U)=(h,f-Ug)$ contains an element with a term of first order. From Theorem~\ref{catlinmultthm}, we know $${\text gen.val}_{\{w_1, \dots, w_{q-1} \}} \:\: D\Big( (\I(U), w_1, \dots, w_{q-1}), x_0 \Big) \leq \prod_{i=q}^n \: D_i (\I(U), x_0),$$ while Equation~\eqref{nval} gives us $$\Delta_n(\I(U), x_0)=D_n(\I(U),x_0) = {\text gen.val}_{\{w_1, \dots, w_{n-1} \}}\:\:  D \Big( (\I(U), w_1, \dots, w_{n-1}), x_0 \Big).$$ The fact that $\I(U)$ contains an element with a term of first order means that generically $D \Big( (\I(U), w_1, \dots, w_{n-1}), x_0 \Big)=1,$ which implies $D_n(\I(U),x_0)=1.$ Since $D_k(\I(U), x_0) \leq D_q (\I(U), x_0)$ for all $k \geq q,$ the result follows.   \qed

\smallskip\noindent {\bf Remark:} The claim $D_n(\I(U),x_0)=1$ appears at the top of p.156 of \cite{catlinsubell}.

\bigskip\noindent Each part of Theorem~\ref{maintheorem} will now be proven in a separate proposition.

\medskip
\newtheorem{prop1}[idealleft]{Proposition}
\begin{prop1}
\label{proposition1}
Let $\Omega$ in $\C^n$ be a domain with $\smooth$
boundary. Let $x_0 \in b \Omega$ be a point on the boundary of
the domain. For any $1 \leq q < n,$ $D_q(b \Omega, x_0) \leq \Delta_q(b \Omega, x_0).$
\end{prop1}

\smallskip\noindent {\bf Proof:} Modulo notational changes, the same proof as for Proposition~\ref{idealleftprop} applies here, but we give it again for completeness. If $\Delta_q(b \Omega, x_0)= \infty,$ then the estimate is obviously true. We thus restrict ourselves to the case when $\Delta_q(b \Omega, x_0)=t < \infty.$ Assume the estimate is false, i.e. $D_q(b \Omega, x_0) >t.$ Since $D_q (M, x_0)$ is defined as the supremum over all $q$-dimensional holomorphic varieties passing through $x_0$ of $\tau (V^q, x_0),$ there exists such a holomorphic variety $V^q$ for which $$\tau (V^q, x_0) = {\text gen.val} \left\{ \tau (V^q \cap S, x_0) \right\} =t' >t,$$ but as we have shown in Proposition~\ref{deltaqgeneric}, $\Delta_q(b \Omega, x_0)$ is generic over the choice of $S,$ so the curves in $V^q \cap S$ enter into the computation of $\Delta_q(b \Omega, x_0).$ Therefore, $$\Delta_q(b \Omega, x_0) \geq t' > t,$$ which is obviously a contradiction. \qed

\medskip
\newtheorem{prop2}[idealleft]{Proposition}
\begin{prop2}
\label{proposition2}
Let $\Omega$ in $\C^n$ be a domain with $\smooth$
boundary. Let $x_0 \in b \Omega$ be a point on the boundary of
the domain. For any $1 \leq q < n,$ if the domain is $q$-positive at $x_0,$ then $$\Delta_q(b \Omega, x_0) \leq 2\left(\frac{ D_q(b \Omega, x_0)}{2} \right)^{n-q}.$$
\end{prop2}

\smallskip\noindent {\bf Proof:} Since the domain is $q$-positive at $x_0,$ as shown at the end of Section~\ref{fincat}, we obtain $$\tau(V^q, x_0) \geq 2 \, \tau(\I(U), V^q) \quad \forall \, U, \: \forall \, V^q.$$ Therefore, $$D_q (b \Omega, x_0)=\sup_{V^q}  \left\{ \tau (V^q, x_0) \right\} \geq 2 \,\sup_{V^q}  \left\{ \tau(\I(U), V^q) \right\}= 2\, D_q (I(U), x_0) \quad \forall \, U.$$ In other words, $$ \frac{D_q (b \Omega, x_0)}{2} \geq D_q (I(U), x_0) \quad \forall \, U$$ and $$\left( \frac{D_q (b \Omega, x_0)}{2}\right)^{n-q} \geq \Big(D_q (I(U), x_0)\Big)^{n-q} \quad \forall \, U.$$ We can now take the supremum on the right over all unitary matrices $U \in \U(N)$ and use Corollary~\ref{catlinmultbdrycorollary} to obtain
\begin{equation*}
\begin{split}   
 \left( \frac{D_q (b \Omega, x_0)}{2}\right)^{n-q}& \geq \sup_{U \in \, \U(N)} \Big(D_q (I(U), x_0)\Big)^{n-q}\\&\geq  \sup_{U \in \, \U(N)} {\text gen. val}_{\{w_1, \dots, w_{q-1} \}}\:\: D\Big( (\I(U), w_1, \dots, w_{q-1}), x_0 \Big).
 \end{split}
 \end{equation*}
 By Proposition~\ref{multgeneric} and Equation~\eqref{qbound},
 \begin{equation*}
 \left( \frac{D_q (b \Omega, x_0)}{2}\right)^{n-q} \geq  \sup_{U \in \, \U(N)}  \:\: \inf_{\{w_1, \dots, w_{q-1} \}} D\Big( (\I(U), w_1, \dots, w_{q-1}), x_0 \Big)\geq \frac{1}{2} \:\: \Delta_q(b \Omega, x_0).
 \end{equation*}
\qed

\medskip\noindent {\bf Proof of Theorem~\ref{maintheorem}:} We put together the results of Proposition~\ref{proposition1} with Proposition~\ref{proposition2} and the fact mentioned in Section~\ref{findap} that a pseudoconvex domain is $q$-positive when $\Delta_q(b \Omega, x_0)<\infty,$ namely Proposition~\ref{pscqproppproposition}. \qed

\medskip Let us now address the issue of sharpness for Theorem~\ref{maintheorem} (ii) via an example. Corollary~\ref{catlinmultbdrycorollary}, a consequence of Catlin's Theorem~\ref{catlinmultthm}, is an essential part of the proof of Proposition~\ref{proposition2} and is responsible for the jump in power that destroys any chance for this type of proof to yield the sharp estimate when $q=1.$ Yet, a result like Corollary~\ref{catlinmultbdrycorollary} cannot be avoided because it relates the only truly well-behaved quantity in the problem $D\Big( (\I(U), w_1, \dots, w_{q-1}), x_0 \Big),$ which is upper semi-continuous with respect to $U$ and generic over the choice of $\{w_1, \dots, w_{q-1} \},$ to $D_q.$ Note that upper semi-continuity is necessary because it allows one to handle the supremum over all unitary matrices produced by polarization, exactly the technique that reduces the problem from the ill-behaved local ring $C^\infty_{x_0}$ to $\oka_{x_0}.$ We shall now show by example that Corollary~\ref{catlinmultbdrycorollary} and Theorem~\ref{catlinmultthm} are both sharp. Let $r= Re\{z_3\} + |z_1^3|^2+|z_2^3|^2$ and $n=3.$ Since $g=0$ in the holomorphic decomposition, no unitary matrices appear in its corresponding holomorphic ideal $\I=(z_3, z_1^3, z_2^3).$ $\Delta_1(\I, 0)=D_1 (\I, 0)=3$ and $D(\I,0)=3^2=9,$ so $D(\I, 0)= \left(D_1(\I,0)\right)^2$ for $q=1.$ To show Theorem~\ref{catlinmultthm} is also sharp, we drop $z_3$ from this ideal, and consider $\I=(z_1^3, z_2^3)$ for $n=2.$ Once again, $\Delta_1(\I, 0)=D_1 (\I, 0)=3,$ and $D(\I,0)=9.$ By Equation~\eqref{nval} above, $D_2(\I, 0)={\text gen.val}_{\{w \}}\:\:  D \Big( (\I, w), 0 \Big)=3,$ so when $q=1$ $$D(\I,0)=9=D_1(\I,0) \cdot D_2(\I,0).$$ For comparison, see Example 2.14.1 on p.624 of \cite{opendangelo} used by D'Angelo to show the inequalities in his Theorem 2.7 were sharp.

We now compute $\Delta_q$ and $D_q$ for $q=1,2,3$ for three ideals of germs of holomorphic functions of three variables. Trivially, $\Delta_1=D_1$ and $\Delta_3 = D_3:$

\smallskip\noindent (a) $\I = (z_1^3+ z_2^3-z_3^3).$ Here $\V(\I)$ is a surface, hence $\Delta_1 (\I, 0)=D_1(\I, 0) = \infty$ and $\Delta_2 (\I, 0)=D_2(\I, 0) = \infty,$ while $\Delta_3 (\I, 0)=D_3(\I, 0)  ={\text gen.val}_{\{w_1, w_2 \}}\:\:  D \Big( (\I, w_1,w_2), 0 \Big)= 3.$

\smallskip\noindent (b) $\I = (z_1^3+ z_2^3-z_3^3,(z_1-z_3)^m),$ where $m \in \N,$ $m>3.$ Here $\V(\I)$ is a curve, hence $\Delta_1 (\I, 0)=D_1(\I, 0) = \infty.$ $\Delta_2 (\I, 0)=D_2(\I, 0) = m>3$ is obtained by using $\V((z_1^3+ z_2^3-z_3^3)),$ the variety corresponding to the ideal $\tilde \I = (z_1^3+ z_2^3-z_3^3).$ $\Delta_3 (\I, 0)=D_3(\I, 0) = {\text gen.val}_{\{w_1, w_2 \}}\:\:  D \Big( (\I, w_1,w_2), 0 \Big)= 3.$

\smallskip\noindent (c) Let $\I$ be any ideal generated by homogeneous polynomials in $z_1,$ $z_2,$ and $z_3$ all of degree $p$ satisfying that $\V(\I)=\{0\}.$ In this case, $\Delta_1 (\I, 0)=D_1(\I, 0) = p,$ $\Delta_2 (\I, 0)= D_2(\I, 0) =p,$ while $\Delta_3 (\I, 0)=D_3(\I, 0) = {\text gen.val}_{\{w_1, w_2 \}}\:\:  D \Big( (\I, w_1,w_2), 0 \Big)=p.$

\medskip\noindent In these examples $\Delta_2 = D_2$ as well, but it is the authors' hope that a future investigation will reveal whether equality holds in general or is merely an artifact here of the difficulty of computing $D_q.$

\bibliographystyle{plain}
\bibliography{CatlinDAngeloType}

\end{document}